\documentclass[english,american]{article}
\usepackage[T1]{fontenc}
\usepackage{textcomp}
\usepackage[utf8]{inputenc}
\usepackage{babel}
\usepackage{booktabs}
\usepackage{varwidth}
\usepackage{tabularx}
\usepackage{amssymb}
\usepackage[letterpaper]{geometry}
\usepackage{setspace}
\usepackage[]
 {hyperref}

\makeatletter

\providecommand{\tabularnewline}{\\}
\newenvironment{cellvarwidth}[1][t]
    {\begin{varwidth}[#1]{\linewidth}}
    {\@finalstrut\@arstrutbox\end{varwidth}}

\@ifundefined{date}{}{\date{}}
\makeatother

\begin{document}
\title{Using Large Language Models to Study Mathematical Practice}
\author{William D'Alessandro}
\maketitle
\begin{abstract}
\noindent The philosophy of mathematical practice (PMP) looks to evidence
from working mathematics to help settle philosophical questions. One
prominent program under the PMP banner is the study of explanation
in mathematics, which aims to understand what sorts of proofs mathematicians
consider explanatory and what role the pursuit of explanation plays
in mathematical practice. In an effort to address worries about cherry-picked
examples and file-drawer problems in PMP, a handful of authors have
recently turned to corpus analysis methods as a promising alternative
to small-scale case studies. But this research has itself been criticized
for attempting to draw philosophical conclusions from counts of ambiguous
terms like `explain'. What’s needed to propel such methods to greater
success, it seems, is a language-processing tool capable of taking
in large amounts of text and reasoning insightfully about the meanings
it contains. The current generation of frontier language models is
already adept at sophisticated text analysis. The main bottleneck,
now on the verge of being overcome, has been the capability for researchers
to feed such systems big-data-sized text corpora (and to receive reliably
correct outputs in return). This paper reports the results from such
a corpus study facilitated by Google's Gemini 2.5 Pro, a model whose
reasoning capabilities, advances in hallucination control and large
context window allow for the accurate analysis of hundreds of pages
of text per query. Based on a sample of 5000 mathematics papers from
arXiv.org, the experiments yielded a dataset of hundreds of useful
annotated examples. Its aim was to gain insight on questions like
the following: How often do mathematicians make claims about explanation
in the relevant sense? Do mathematicians' explanatory practices vary
in any noticeable way by subject matter? Which philosophical theories
of explanation are most consistent with a large body of non-cherry-picked
examples? How might philosophers make further use of AI tools to gain
insights from large datasets of this kind? As the first PMP study
making extensive use of LLM methods, it also seeks to begin a conversation
about these methods as research tools in practice-oriented philosophy
and to evaluate the strengths and weaknesses of current models for
such work.
\end{abstract}
{[}Draft of an invited chapter for \emph{Mathematicians at Work: Empirically
Informed Philosophy of Mathematics}, ed. by José Perez Escobar, Deborah
Kant, Deniz Sarikaya and Mira Sarikaya, \emph{Synthese Library} series{]}

\section{Introduction}

The philosophy of mathematical practice (PMP) looks to working mathematics
as a primary source of guidance and evidence for philosophical theorizing.
The study of explanation in mathematics is a case in point. It was
set in motion by the observation that mathematicians often distinguish
explanatory from non-explanatory proofs (\cite{Steiner 1978}), and
work on the topic in a PMP vein continues to appeal to the mathematical
literature to supply motivating examples and adjudicate between rival
theories. 

The core methodological problem of PMP is that mathematical practice
can only serve in this advisory role to the extent that we understand
and agree about the details of the relevant practices. How best to
gain this understanding isn't a priori obvious. Indeed, the question
of proper methods for PMP has recently been a topic of debate (\cite{Rittberg 2019},
\cite{D=002019Alessandro 2025}).

The default approach throughout most of the subject's history---following
in the footsteps of similar practice-oriented paradigms in philosophy
of science and elsewhere---has granted pride of place to case studies
focused on small numbers of texts or issues. At their carefully-drawn
best, such studies shine revealing light on important mathematical
phenomena. Yet the case study method has limits. For one, the detailed
scholarship required means that high-quality studies are inherently
constrained in volume, so they can at best provide only a narrow window
on a practice of interest. Moreover, selection effects presumably
influence case study choice in various ways: especially famous, striking,
controversial or novel cases can be expected to garner disproportionate
interest from researchers and publishing gatekeepers, and these biases
may lead to a misleading picture of practice as a whole. Finally,
researchers may be tempted to strategically highlight examples that
support their favored theories or otherwise suit their individual
purposes.

Motivated by such worries, a number of philosophers have called for
a reevaluation of the case study as PMP's default methodology. Alongside
it---or perhaps even in its place---some hope to elevate a more
objective, systematic, data-driven paradigm, less vulnerable to nagging
worries about cherry-picking and file-drawer problems. This reformist
vision might include any of several elements. Some authors, for instance,
have embraced the tools of quantitative social science as means to
study mathematicians' judgments experimentally (\cite{Inglis =000026 Aberdein 2015},
\cite{Inglis =000026 Aberdein 2016}, \cite{Mej=0000EDa Ramos et al. 2021},
\cite{Sa et al. 2024}). Others have applied big-data methods to mathematical
corpora in pursuit of insights about practices in aggregate (\cite{Johansen et al. 2018},
\cite{Mej=0000EDa Ramos et al. 2019}, \cite{Pease et al. 2019},
\cite{Mizrahi 2020}).\footnote{This isn't to say that traditional case studies and data-driven methods
are the only live approaches to PMP. On the contrary, \cite{Rittberg 2019}
documents seven distinct PMP research frameworks, including interviews,
cognitive-scientific methods and collaborations between mathematicians
and philosophers (13).}

This latter tradition supplies the motivation and backdrop for the
present study. In previous work, large-scale corpus analysis in the
PMP setting has been based primarily on counting occurrences of key
terms and doing statistics on such counts. But this approach has significant
weaknesses. \foreignlanguage{english}{A central issue is that ``indicator
phrases'' related to topics of PMP interest are highly ambiguous
and context-sensitive; this is the ``word-concept problem'' for
corpus methods discussed by \cite{Chartrand 2022}. Consider explanation,
for instance. Terms like `explains' and its cognates sometimes indicate
genuine explanations (in the sense of reasons why a mathematical fact
is true), but they more often indicate other notions: explaining what
a term means, how to perform a calculation, why a notation was chosen,
how to interpret a theorem, and so on. In short, corpus work based
on simple keyword counts ``interprets the use of words like `explain'...
as indicating that the text excerpt that contains this segment is
an explanation. But obviously people do not label their discourse
as they are engaging in it'' (\cite{Chartrand 2022}, 8). Conversely,
there are many ways to raise explanatory concerns in alternative language.
And of course the same can be said for other philosophically interesting
notions: beauty, simplicity, depth, understanding, naturalness and
so on. }

As a result, corpus methods' potential to confidently draw substantive
PMP conclusions has thus far been limited. While the basic word-count
paradigm can undoubtedly be improved upon---\foreignlanguage{english}{e.g.,
researchers can choose focal terms with greater care, tailor hypotheses
more narrowly, or individually check whether each data-mined item
falls under the concept of interest---the challenges along this path
are formidable, and the knowledge to be gained seems quite constrained
even under ideal conditions. }

\selectlanguage{english}%
As \cite{D=002019Alessandro 2025} suggests, the use of large language
models (LLMs) for automated text analysis may represent an important
path forward for corpus methods. The key advantage of LLMs is their
proficiency at grasping, reasoning about and reporting on the meanings
of complex texts.\footnote{\selectlanguage{american}%
I take no stand here on whether LLMs have genuine mental states or
cognitive capabilities. In saying that they grasp and reason about
the meanings of texts, I mean that their outputs are comparable to
those of sophisticated humans engaged in these activities.\selectlanguage{english}%
} In particular, given an appropriate prompt, these models can often
readily identify whether a concept of interest appears in a particular
discussion, regardless of the specific language in which the concept
is invoked. This functionality seems to promise much improved accuracy,
flexibility and analytical depth to researchers in the corpus tradition.

\selectlanguage{american}%
It's only quite recently, however, that such hopes have begun to
seem fully realizable. Early-generation LLMs suffered from several
problems that limited their research utility. One longstanding issue
is \emph{context window length}, the maximum input size a model can
accept.\emph{ }The launch versions of ChatGPT and GPT-4 featured context
lengths of under $10,000$ tokens (\textasciitilde$7500$ words),
less than the length of many individual research papers and not nearly
enough for large-scale corpus work. Context windows have generally
grown since, but only in recent months have frontier models begun
to boast usefully large windows of a million tokens or more (\textasciitilde$750,000$
words, or \textasciitilde 2000 pages). The performance of leading
models on longer inputs has also improved significantly:\footnote{As measured, for instance, by benchmarks like LongBench; see https://longbench2.github.io/\#leaderboard.}
accurately retrieving specific facts embedded in large prompts (in
so-called needle-in-a-haystack tasks) was virtually impossible for
early models but now appears to be a solved problem.\footnote{For instance, Google reports that its Gemini 2.5 Pro achieves $94.5\%$
on the MRCR long-context benchmark (https://storage.googleapis.com/model-cards/documents/gemini-2.5-pro-preview.pdf).
See https://github.com/gkamradt/LLMTest\_NeedleInAHaystack for some
results on earlier models.} Finally, LLMs' mathematical capabilities---and complex cognition
more generally---were greatly enhanced with the arrival of reasoning
models, beginning with OpenAI's o1 in late 2024. 

This suite of developments has made LLM-assisted corpus work begin
to look feasible. Current large-context models can accept long prompts
featuring detailed instructions and at least dozens of papers to analyze.
They can moreover be counted on to correctly report and intelligently
reason about the contents of these texts. With the ability to automate
queries via API interfaces, researchers are in a position to analyze
large batches of material in this way with a minimum of human labor. 

The present project set out to test the waters. It addresses a number
of foundational questions for future PMP research. The first is logistical:
Given a big-data-sized corpus of mathematical texts (comprising several
thousand papers, say), what sort of pipeline can be devised to efficiently
feed the corpus to a chosen LLM in batches while receiving running
reports on results for each batch? The second is a question of prompt
engineering: Is it possible to devise instructions to the model which
will yield useful and accurate reports, identifying as many genuine
occurrences of the target concept as possible while minimizing false
positives? The third question is epistemic: Given that the method
can be made to work, what can a dataset obtained in this way tell
us about questions of PMP interest? The fourth is about LLMs' nascent
philosophical capabilities: can these models themselves, equipped
with the data from a large-scale corpus analysis, help settle debates
between rival theories, or even suggest entirely new lines of inquiry?

As the discussion below will show, optimism about LLM methods in PMP
(and elsewhere) is in many respects warranted. Current models are
capable of doing sophisticated, accurate and interesting corpus work
on a scale that would be impossible by other means. Their ability
to intervene in philosophical controversies is less clear, but preliminary
results suggest potential roles for AI tools even here. 

\section{Methods}

The corpus for the study consisted of $5000$ PDFs from the Mathematics
category on arXiv.org (a random sample from the full arXiv Mathematics
repository of around\textbf{ $80,000$} papers, downloaded from Google
Cloud Storage using the $\mathtt{gsutil}$ tool). arXiv preprints
are a key medium for research dissemination in contemporary mathematics,
as in other hard sciences, which often precede or complement traditional
publishing. Like the contents of the arXiv as a whole, the papers
in the corpus represented a very broad spectrum of mathematical subjects,
styles and authors. I view this as a virtue of a corpus aiming to
impart lessons about mathematical practice generally (as opposed to
reflecting, say, only the type of research published in a handful
of top journals). Note that, since the arXiv has featured a Mathematics
category only since the mid-1990s, one limitation of the corpus is
that it represents a fairly narrow and recent time-slice of mathematical
activity. 

The corpus analysis primarily made use of Google's Gemini 2.5 Pro
Experimental and 2.5 Pro Preview, released on March 25 2025 and June
5 2025 respectively. Gemini 2.5 Pro was by many standards the best-performing
publicly available model at the time of release (and remains so at
the time of writing), with a context window of one million tokens
and integrated chain-of-thought reasoning. Initial tests were done
with Gemini 1.5 Pro on account of its still larger two-million-token
context, but I found this earlier model prone to giving inaccurate
and hallucinated responses on long inputs. 2.5 Pro exhibited no signs
of this problem; the model's references to source materials proved
accurate on each of the many spot checks I performed over the course
of the project.

My goal in prompting the model was to give abundant context and examples
for the target concept of explanation without introducing controversial
assumptions unduly favoring my own (or some other narrow set of) philosophical
views.\footnote{The point has been made that not all explanations in mathematics take
the form of traditional proofs (\cite{Lange 2017}, \cite{D=002019Alessandro 2020}).
Some have moreover argued that not all explanations are best understood
as giving the reason why a particular fact holds (\cite{Lehet 2021},
\cite{D=002019Alessandro =000026 Lehet 2024}). The prompt used for
the present study nevertheless focuses on traditional cases of explanatory
proof in order to keep complexity and controversy to a minimum.} To achieve these aims, I used a prompt consisting mostly of a long
(\textasciitilde 5000-word) excerpt from \cite{Mancosu et al. 2023},
the \emph{Stanford Encyclopedia of Philosophy} article on explanation
in mathematics by Paolo Mancosu, Francesca Poggiolesi and Chris Pincock.
This excerpt---consisting of sections 2, 2.1 and 2.2 of the article---includes
a historical overview of explanatory concerns in mathematics and an
overview of major contemporary views on the topic, interspersed with
plentiful examples. The excerpt text was prefaced with a short summary
introduction and some instructions to the model on how to perform
the analysis and structure its outputs. The complete prompt was as
follows:
\begin{quote}
``You are a skilled and insightful research assistant on a text analysis
project. You’re knowledgeable about research-level mathematics. You’re
also well-versed in philosophy, attuned to subtle conceptual distinctions
and proficient at close reading. You take care to properly understand
the meaning of texts and always carefully state and justify your interpretations. 

Your task is to analyze mathematical research papers to identify instances
and discussions of mathematical explanation within pure mathematics.
\emph{Mathematical explanation} refers to cases where mathematicians
go beyond merely establishing the truth of results and seek to provide
insight into why mathematical facts hold. 

Some examples of phenomena to watch for include the following:
\begin{itemize}
\item Explanatory vs. non-explanatory proofs: Look for discussions where
mathematicians distinguish between proofs that merely establish that
a result is true versus those that explain why it is true, or provide
explanatory understanding of its truth. Document when mathematicians
contrast different approaches to the same theorem, especially when
one approach is favored for explanatory reasons.
\item Reproof efforts: Identify cases where theorems are proved multiple
times, with the search for additional proofs justified on explanatory
grounds. 
\item Explanation-related terminology: Search for terms like `explain',
`reason why', `understanding why', `insight', `intuition'
or `deep reason' in relation to mathematical proofs or results.
\item Explanatory failure terminology: Look for critiques of proofs as unexplanatory,
for instance using descriptions like `brute-force', `not illuminating'
or `unintuitive'.
\end{itemize}
Some types of potentially explanatory proofs: 
\begin{itemize}
\item Visual or diagrammatic proofs that provide insight
\item Proofs that use characteristic properties of the entities involved
\item Proofs that exploit symmetry, unity, or simplicity
\item Proofs that reveal underlying structures or mechanisms
\item Proofs that unify different cases.
\end{itemize}
Pay particular attention to the following important instructions. 
\begin{itemize}
\item When analyzing the attached mathematics papers, provide direct quotes
from mathematicians discussing explanatory value, along with context
about the mathematical content they're addressing. Focus on how mathematicians
themselves articulate the difference between knowing that a result
is true and understanding why it is true. 
\item Give the filename, title, author and page number for each reference
you mention. Do not hallucinate content, confuse sources or make up
quotes.
\item In addition to explicit explanatory success claims, you should note
and report claims about explanatory failures, desires or hopes for
explanation, possible routes to explanation, and the like. But be
sure to only include relatively clear examples rather than vaguely
or tangentially related points.
\item Don't summarize the content of each individual paper. Only report
on findings that are potentially relevant to the prompt. If a paper
contains no such examples, don't mention it in your analysis. Overall,
be efficient, and prioritize reporting every relevant example over
analyzing each case at length. 
\end{itemize}
The remainder of the prompt is an excerpt from Mancosu et al.'s survey
article, `Mathematical Explanation'. It provides further context
on the target concept(s) of explanation. Don't discuss this article
or the examples it contains in your analysis of the texts. BEGINNING
OF MANCOSU ET AL. ARTICLE EXCERPT: {[}The remainder of the prompt
consisted of sections 2, 2.1 and 2.2 of \cite{Mancosu et al. 2023}.{]}''
\end{quote}
Even with Gemini 2.5's large context window, processing the whole
corpus of 5000 PDFs with the above prompt required many individual
runs, each requiring several minutes of thinking time and output text
generation. In order to automate this process, I used a Python script
produced by Gemini itself. Impressively, this (rather complex) code
worked exactly as needed, allowing me to set up the automation pipeline
with minimal difficulty in spite of my scant personal coding experience.
After each run, Gemini saved its analysis of the relevant batch of
PDFs to a text file. Processing the entire corpus required 200 runs
(of 25 papers each) and approximately 24 hours.

This process yielded 200 files containing Gemini's analyses of the
corpus texts. Each analysis took the form of a list of candidate examples
of explanation or explanatory concerns appearing in one batch of texts.
The reported example items adhered to the format {[}filename, title,
authors, summary of example, quote, Gemini's commentary{]}. A randomly
selected example follows. (Note that here and throughout the paper,
text produced by Gemini is rendered in sans-serif type.\footnote{I've lightly copyedited these items for uniformity and clarity, as
Gemini's outputs across different runs vary slightly in formatting.
The ``example'' and ``context'' text in each case is verbatim
from Gemini.})
\begin{itemize}
\item \textsf{\textbf{Title:}}\textsf{ From Mennicke Symbols to Euler Class
Groups (\cite{van der Kallen 2000})}
\item \textsf{\textbf{Example: }}\textsf{An analogy with topology is cited
as providing an explanation for an algebraic structure.}
\item \textsf{\textbf{Quote: }}\textsf{``Let us now take $A$ to be the
Banach algebra of continuous real valued functions on some finite
$d$-dimensional CW complex $X$. Then one knows that for $n\ge3$
the orbit set $\mathrm{Um}_{n}\left(A\right)/E_{n}\left(A\right)$
is in bijective correspondence with the set $\left[X,\mathbb{R}^{n}-0\right]$
of homotopy classes of maps from $X$ to $\mathbb{R}^{n}-0=\mathrm{Um}_{n}\left(\mathbb{R}\right)$.
This gives a topological explanation why for $2\le d\le2n-4$ one
has a group structure on $\mathrm{Um}_{n}\left(A\right)/E_{n}\left(A\right)$''
(p. 10).}
\item \textsf{\textbf{Context:}}\textsf{ Discussing orbit sets over Banach
algebras $A=C(X)$. The author explicitly labels the connection to
homotopy theory as a ``topological explanation'' for the existence
of a group structure.}
\end{itemize}
A second example:
\begin{itemize}
\item \textsf{\textbf{Title:}}\textsf{ Formulas of Verlinde Type for Non-Simply
Connected Groups (\cite{Alekseev et al. 2000})}
\item \textsf{\textbf{Example: }}\textsf{Page 1 (Introduction): States the
motivation is to apply the fixed point formula from the companion
paper to understand Verlinde's formula for geometric quantization
of moduli spaces, connecting index theory on loop group spaces to
formulas arising in conformal field theory and algebraic geometry.}
\item \textsf{\textbf{Quote: }}\textsf{``In this paper we give applications
of the fixed point formula proved in the companion paper. Our original
motivation was to understand a formula of E. Verlinde for the geometric
quantization of the moduli space of flat connections on a Riemann
surface. In particular A. Szenes suggested to us that the Verlinde
formula should follow from an equivariant index theorem, much as the
Weyl or Steinberg formulas can be interpreted as fixed point formulas
for flag varieties'' (p. 1).}
\item \textsf{\textbf{Context:}}\textsf{ This explicitly frames the work
as seeking an explanation (``understand a formula of E. Verlinde'')
by deriving it from a more general principle (equivariant index theorem
/ fixed point formula), thus providing deeper insight into the Verlinde
formula's origins and connections, explaining why it holds.}
\end{itemize}
A third: 
\begin{itemize}
\item \textsf{\textbf{Title:}}\textsf{ Modular Data: The Algebraic Combinatorics
of Conformal Field Theory (\cite{Gannon 2002})}
\item \textsf{\textbf{Finding: }}\textsf{Discussion of seeking underlying
reasons for observed patterns.}
\item \textsf{\textbf{Quote: }}\textsf{``Patterns such as A-D-E are usually
explained by identifying an underlying combinatorial fact which is
responsible for its various incarnations. The A-D-E combinatorial
fact is probably the classification of symmetric matrices over $\mathbb{Z}_{\geq}$
... Perhaps the only A-D-E classification which still resists this
‘explanation’ is that of $A_{1}^{\left(1\right)}$ modular invariants''
(p. 29).}
\item \textsf{\textbf{Context:}}\textsf{ The author discusses the recurrence
of A-D-E classification schemes in various mathematical contexts and
notes that these patterns are typically {*}explained{*} by finding
a common underlying combinatorial structure. The quote highlights
the search for such an explanation for the $A_{1}^{\left(1\right)}$
modular invariants, noting it as a current explanatory gap.}
\end{itemize}
A fourth:
\begin{itemize}
\item \textsf{\textbf{Title:}}\textsf{ Composition Sum Identities Related
to the Distribution of Coordinate Values in a Discrete Simplex (\cite{Milson 2000})}
\item \textsf{\textbf{Finding: }}\textsf{Explaining the reason behind a
mathematical property (exact solvability) by relating it to a known
structure or equivalence.}
\item \textsf{\textbf{Quote: }}\textsf{``Interesting composition sum identities
will appear in the present context when we consider exactly-solvable
differential equations. We present three such examples below, and
discuss the enumerative interpretations in the next section. In each
case the exact solvability comes about because the equation is gauge-equivalent
to either the hypergeometric, or the confluent hypergeometric equation''
(p. 8).}
\item \textsf{\textbf{Context:}}\textsf{ Introducing three examples of second-order
differential equations whose series solutions lead to composition
sum identities (Propositions 4.2, 4.3, 4.4). The author explains {*}why{*}
these specific equations are exactly solvable, attributing it to their
gauge-equivalence to standard, well-understood hypergeometric equations.}
\end{itemize}
And a fifth:
\begin{itemize}
\item \textsf{\textbf{Title:}}\textsf{ Combinatorial Formulas for Products
of Thom Classes (\cite{Guillemin =000026 Zara 2000})}
\item \textsf{\textbf{Finding: }}\textsf{The authors are discussing the
organization of the paper and highlighting a particularly interesting
aspect of their formula (1.11) for Thom classes in equivariant cohomology.}
\item \textsf{\textbf{Quote: }}\textsf{``In Section 5 we will attempt to
demystify what is perhaps the most puzzling feature of the formula
(1.11), the fact that all the summands are rational functions (elements
of the quotient field, $Q\left(\mathfrak{g}^{*}\right)$), whereas
the sum itself is a polynomial. This indicates that a lot of mysterious
cancellations are occurring in this summation; and we will show }\textsf{\emph{how}}\textsf{
these cancellations occur in a few simple but enlightening examples''
(p. 6).}
\item \textsf{\textbf{Context:}}\textsf{ The terms ``demystify'', ``puzzling
feature'', ``mysterious cancellations'', and ``enlightening examples''
strongly indicate a concern for explanation. The authors acknowledge
that the formula, while correct, has a feature that lacks immediate
understanding (why rational functions sum to a polynomial). They explicitly
aim to provide insight into the }\textsf{\emph{reason why}}\textsf{
this happens by analyzing the cancellation mechanism in simple cases,
moving beyond just knowing }\textsf{\emph{that}}\textsf{ the formula
yields a polynomial to understanding }\textsf{\emph{how/why}}\textsf{
it does. This aligns with seeking explanatory understanding and revealing
underlying mechanisms. }
\end{itemize}
Gemini's 200 output files delivered a very large number of putative
examples with an enormous total length of several hundred thousand
words. Many of these examples were excellent. In particular, one can
see from the above that Gemini succeeded at identifying references
to explanation phrased in a variety of ways, bearing out the promise
of LLM tools to overcome the word-concept problem for traditional
corpus methods. 

It was clear on inspection, however, that Gemini's first pass had
erred considerably on the side of generosity, including many low-quality
cases only marginally or speculatively related to explanation. For
instance, the model often flagged passages which merely mentioned
links between theories or problems, noted generalizations, contained
generic positive evaluations (``well-known'', ``instructive'',
``nice'', ``precise'') and the like, tendentiously asserting the
presence of explanatory concerns in these passages for reasons not
supported by the texts themselves.

In order to improve the quality of the dataset for downstream analysis,
I aimed to rerun Gemini's reports with a new prompt achieving extensive
filtering of weaker items. Initial attempts to do this were unsuccessful;
in almost all cases the model would simply reproduce the original
analyses with no filtering. I achieved better results by asking Gemini
to compose a prompt to elicit the desired behavior. It produced the
following:
\begin{quote}
\textsf{I'm providing you with a set of text files from a previous
Gemini analysis. Your sole and critical task now is to act as a strict
quality filter for these files.}

\textsf{The original analysis (which you should not try to replicate)
aimed to find examples of mathematical explanation, defined as cases
where mathematicians go beyond merely establishing truth and seek
to provide insight into why mathematical facts hold. This includes
distinguishing explanatory vs. non-explanatory proofs, reproofs for
explanatory reasons, explicit discussion of `explain', `reason
why', `understanding why', `insight into why', etc.}

\textsf{Your current goal: Read through the provided snippets from
the previous Gemini output. Aggressively filter these snippets. Select
only the most clear, unambiguous, and explicit references to mathematical
explanation or explanatory concerns. You MUST exclude at least 50-60\%
of the original examples. Many items in the input files are only marginally,
speculatively, or not at all related to explanatory concerns as defined.}

\textsf{Your guiding principle for exclusion: If there's a plausible
interpretation of the example where the authors aren't invoking specifically
explanatory concerns, then OMIT THE EXAMPLE. If in any doubt, leave
it out. Only include an example if the most direct and unavoidable
interpretation requires concepts directly related to understanding
why a mathematical statement is true.}

\textsf{Retain original content and formatting for the examples you
do select. Do not include lists of files processed. Each time you
start listing examples from a new \textquotedbl batch\_n\_output.txt\textquotedbl{}
file, name that file.}

\textsf{Key Distinctions for Filtering (Apply these strictly):}

\textsf{KEEP: }
\begin{itemize}
\item \textsf{Explicit mentions of ``explaining why'', ``reason for'',
``understanding the mechanism/structure behind'', ``deeper reason'',
distinguishing proofs that ``show why'' versus those that merely
``show that'', reproving for ``better understanding''.}
\end{itemize}
\textsf{DISCARD: }
\begin{itemize}
\item \textsf{General claims of ``understanding'', ``insight'', ``intuition'',
``clarity'', ``simplicity'' if not clearly tied to understanding
why a result is true (e.g., understanding how to do something is not
the target).}
\item \textsf{Mentions of analogies, reduction strategies, comparisons with
alternatives, or subsuming cases under a general principle, UNLESS
these are explicitly linked to providing an explanation for why something
is true. Many non-explanatory proofs use these.}
\item \textsf{Statements about a concept being ``instructive'', ``useful
for construction'', ``reducing a problem'', or ``playing a crucial
role'' if they don't also clearly articulate an explanatory benefit
related to why the mathematics is the way it is.}
\end{itemize}
\textsf{Here are some reference examples. Use these to calibrate your
strictness.}

\textsf{Clear instances to REPORT (these exemplify the high bar): }
\begin{itemize}
\item \textsf{``We want to understand the homomorphism and the subgroup
in other ways. First of all, we have in {[}vdK2{]} a more mysterious
construction of the group structure, by generators and relations.
This more mysterious construction works as soon as $2\le d\le2n-4$.
So one would expect to construct a homomorphism... by checking that
the relations hold in $E(A)$. When trying to do this, one finds that
the defining relations used in {[}vdK2{]} are not so obvious in $E(A)$.''
(Context: explicit desire to understand a ``mysterious'' construction
better).}
\item \textsf{``The reason for this fact... at first glance quite surprising,
is of a topological nature...'' (Context: giving a reason for a surprising
fact).}
\item \textsf{``I hope to exhibit the common mechanism underlying...''
(Context: revealing common mechanisms is often explanatory).}
\item \textsf{``We do not know of an argument which provides a unified
explanation...'' (Context: explicit statement of explanatory failure/gap).}
\end{itemize}
\textsf{Unclear instances to OMIT (these are examples of what to discard):}
\begin{itemize}
\item \textsf{``The idea of introducing associative triples is to reduce
the problem... and create a tool for the induction process.'' (No
clear link to explanatory understanding).}
\item \textsf{``It is instructive to reformulate the data...'' (Instructive
doesn't inherently mean explanatory).}
\item \textsf{``The macroscopic effect... can be described like this...
We give a precise definition... whose association... is well-known.''
(Description/association, not necessarily explanation).}
\item \textsf{``A crucial role throughout the paper is played by the decay
behavior...'' (``Crucial role'' is not necessarily an explanatory
role).}
\end{itemize}
\textsf{Remember, your primary function here is selective filtering
based on very strict criteria for what counts as an explicit and unambiguous
discussion of mathematical explanation. Be ruthless in your exclusions.}
\end{quote}
This prompt caused Gemini to remove a large proportion of the examples
from its original outputs. The remaining \textasciitilde 1250 examples
from \textasciitilde 735 distinct sources were of generally higher
quality, and inspection of the model's chain of thought\footnote{Via the ``view model thoughts'' functionality in Google AI Studio,
currently labeled ``experimental''.} suggested that it understood the filtering task correctly. For example,
the model chose to omit a case mentioning the ``remarkable properties''
of residue complexes on the grounds that ``\textsf{`remarkable'
means `powerful' or `useful', not necessarily `explanatory of
why the theory works'.}'' Meanwhile, Gemini chose to retain a case
from \cite{Massey 2001} in which an earlier proof was criticized
as ``too magical'' (``in no way does one `see' the contributions
from the critical points of a small perturbation of $f$'', 2). The
model's reasoning was that this reference ``\textsf{contrast{[}s{]}
an effective but non-explanatory proof style... with an implied desire
for a proof that provides more insight or intuition into the mechanism}''.
The full dataset is available as a (500+-page) Word document at \href{https://www.dropbox.com/scl/fi/djh03oqmfr1mz16msn2ez/Gemini-filtered-outputs.docx?rlkey=0nz5xkwsq8rd8wi1uhexdd3z3&st=74anpvta&dl=0}{this URL}. 

Even after substantial filtering, the dataset remains imperfect. While
the model succeeded at culling a large number of low-quality cases
for appropriate reasons reflecting a good understanding of the task,
it appears unwilling to change its mind about other equally weak cases.
Repeated applications of the filtering prompt to the previously-filtered
dataset led to no noticeable improvements here. The reason for this
behavior is unclear to me.

Moreover and less surprisingly, many cases in the filtered dataset
represent a gray area where something like an explanatory concern
is plausibly present but the language is vague or ambiguous enough
to permit doubt. (Example: \cite{Allcock et al. 2000} described its
methods as ``giv{[}ing{]} a different insight into the structure
of the fundamental group of the space of smooth cubic surfaces''
(3). Gemini retained this case on the grounds that it ``\textsf{addresses
the goal of seeking understanding and insight beyond just determining
the group's structure}''. But it's in my view unclear whether the
insight in question represents an explanatory achievement, as opposed
to mere descriptive knowledge or objectual understanding.) Whether
such cases should be included or excluded from a PMP study of explanation
is in part a matter of taste and in part a question of the goals of
one's specific project. Of course, given that such tastes and goals
vary among human experts, even the best automated classification will
inevitably be somewhat controversial.

After filtering I estimate that Gemini's dataset consists of roughly
20\% high-quality cases (in which substantive explanatory concerns
are unmistakably on display), roughly 20\% low-quality cases (in which
there's little clear sign of such concerns), and roughly 60\% borderline
cases which one might or might not wish to count for various purposes.
Given the large size of the dataset, these percentages represent an
impressive number of useful or potentially useful examples---perhaps
around 250 high-quality and 750 medium-quality cases. Even the smaller
set of top-shelf examples likely rivals or exceeds the total number
of case studies in the literature to date.

\section{LLMs as analytical tools}

A further advantage of LLM-assisted corpus work over traditional approaches
is that, in addition to their utility for amassing large repositories
of examples, language models can also be queried for analysis of these
repositories. These analyses might involve relatively straightforward
assessments of objective features, or they might veer into interpretative
philosophical territory. I report on a few of these applications below:
using Gemini to quantify explanatory practices by subfield (§3.1)
as well as the overall prevalence of explanation-talk (§3.2), querying
Gemini's dataset to find targeted examples of interest (§3.3), and
enlisting Gemini as a philosophical assistant in various capacities
(§3.4).

\subsection{Explanatory practices by subfield}

To start with, PMP researchers might want to know which \emph{subfields
of mathematics} are most or least likely to feature appeals to explanation.
Such high-level features of mathematical practice are both potentially
very interesting and difficult to glean from individual case studies.
According to Gemini's analysis, the predominating fields in the dataset
include geometry, algebra, topology and and combinatorics. As Gemini
explains at length\footnote{Here and below, I've converted Gemini's filename references to standard
citations; the output is otherwise unchanged.}:
\begin{quote}
\textsf{These fields share a common feature: their objects of study
are incredibly rich and can be viewed from multiple perspectives (e.g.,
algebraic, geometric, combinatorial). This multiplicity of perspectives
naturally fosters explanatory questions. The appeals tend to take
the following forms:}
\begin{itemize}
\item \textsf{\textbf{Unification and Revealing Deeper Structures:}}\textsf{
This is the most common form of explanatory claim. A result is deemed
explanatory when it shows that two or more seemingly disparate concepts
are, in fact, instances of a single, deeper underlying structure.}
\begin{itemize}
\item \textsf{\textbf{Example (Algebraic Geometry/Number Theory): }}\textsf{The
modularity conjecture is presented as explaining the analytic properties
of Artin $L$-functions by revealing their origin in the deeper, more
structured world of automorphic forms (\cite{Ramakrishnan 2001}).}
\item \textsf{\textbf{Example (Representation Theory):}}\textsf{ The work
on double affine Hecke algebras is explicitly framed as a ``unification''
of different transforms, which restores desirable properties (like
self-duality) that were lost in previous generalizations. This new
framework explains the properties of the transforms by revealing their
common origin (\cite{Cherednik 2000}).}
\end{itemize}
\item \textsf{\textbf{Contrasting and Justifying Proof Styles: }}\textsf{Mathematicians
in these areas frequently contrast different proof methods, valuing
one over another for its explanatory power. A proof is often praised
for being more ``conceptual'', ``geometric'', ``elementary'',
or ``motivated'', and critiqued for being ``computational'', ``a
case-by-case study'', ``magical'', or ``unenlightening''.}
\begin{itemize}
\item \textsf{\textbf{Example (Topology):}}\textsf{ A new, ``far less computational''
proof for properties of the Birman-Ko-Lee monoid is justified because
its ``geometric interpretation'' provides more insight than the
original, more calculation-heavy one (\cite{Bessis et al. 2000}).}
\item \textsf{\textbf{Example (Algebraic Geometry):}}\textsf{ A proof of
the existence of Whitney Stratifications by Wall is described as \textquotedbl the
most illuminating\textquotedbl{} because it \textquotedbl simplifies\textquotedbl{}
prior proofs, contrasting it with Hironaka's \textquotedbl elegant\textquotedbl{}
but less accessible proof that requires heavy background in resolution
of singularities (\cite{Kaloshin 2000}).}
\end{itemize}
\item \textsf{\textbf{Providing a ``Reason Why'' via a Mechanism: }}\textsf{The
goal is often to move beyond verifying a fact to explaining the mechanism
that brings it about.}
\begin{itemize}
\item \textsf{\textbf{Example (Topology/Symplectic Geometry):}}\textsf{
The change in the homotopy type of a symplectomorphism group is explained
not just by observing it, but by identifying the mechanism: the appearance
of a new stratum in a moduli space, which generates a ``fragile''
but non-trivial homotopy class that subsequently vanishes (\cite{McDuff 2004}).}
\item \textsf{\textbf{Example (Algebraic Combinatorics):}}\textsf{ The search
for bijective proofs is a quintessential example. A generating function
identity may prove that two sets have the same size, but a bijection
explains this equality by providing a direct structural correspondence
between the elements of the sets (\cite{Bressoud 2000}).}
\end{itemize}
\item \textsf{\textbf{Explicitly Stating a Lack of Explanation: }}\textsf{These
fields are rife with authors explicitly pointing out explanatory gaps,
using phrases like ``it is unclear why'', ``calls for an intuitive
explanation'', ``this seems mysterious'', or ``we have no `reason'
for it, except calculation''.}
\begin{itemize}
\item \textsf{\textbf{Example (Algebraic Geometry): }}\textsf{The authors
of \cite{Dolgachev =000026 Keum 2001} state, ``The reason why the
beautiful combinatorics of the Leech lattice plays a role... is still
unclear to us. We hope that the classification of all $\mathrm{K}3$-surfaces...
will shed more light to this question''.}
\end{itemize}
\end{itemize}
\end{quote}
Meanwhile, according to Gemini, the areas in the dataset with fewer
appeals to explanation include core subfields of analysis (measure
theory, classical PDEs, functional analysis), number theory and logic.
Interestingly, Gemini suggests that when explanatory concerns do arise
in these areas, they often differ in character from the types of case
just mentioned:
\begin{itemize}
\item \textsf{\textbf{Justification of Method:}}\textsf{ The explanation
often centers on why a particular method, definition, or function
is the right one to use. The goal is to understand the tool itself,
rather than a property of a separate object.}
\begin{itemize}
\item \textsf{\textbf{Example (Analysis):}}\textsf{ The choice of a specific
``gauge function'' in Henstock-McShane integration is described
as ``rather mysterious'' in its standard presentation. The authors
then claim to provide understanding by showing ``how the properties
of Lebesgue points can be used to determine the gauge'', thus explaining
its origin and role (\cite{Loeb =000026 Talvila 2001}).}
\item \textsf{\textbf{Example (PDE): }}\textsf{In proving regularity for
a wave equation, the author explains the failure of standard iterative
methods (due to a logarithmic divergence) and motivates the necessity
of a new ``adapted frame'' method that is specifically designed
to handle this obstacle (\cite{Tao 2000}). The explanation is about
the success or failure of the proof technology.}
\end{itemize}
\item \textsf{\textbf{Explaining via Reduction to a Fundamental Principle:}}\textsf{
Rather than unifying vast fields, the explanation often consists of
showing how a complex result is a direct consequence of a well-known,
foundational principle within that field.}
\begin{itemize}
\item \textsf{\textbf{Example (Analysis):}}\textsf{ The smoothness and decay
properties of a limiting density function are explained by showing
that its Fourier transform (the characteristic function) has corresponding
properties, and then invoking the standard, powerful theorems of Fourier
analysis (\cite{Fill =000026 Janson 2000}). The explanation is a
clear, structural deduction within an established framework.}
\end{itemize}
\end{itemize}
This dichotomy strikes me as a novel observation deserving of further
attention. I'm aware of little philosophical work addressing the
success or failure of problem-solving techniques or proof technology
as a locus of explanatory concerns.

For the purpose of determining which subfields punch above or below
their weight in explanation-talk, it would be useful to know how the
distribution of papers by area in the filtered dataset compares to
the distribution in the original corpus. To investigate this issue
I ran a script (also produced by Gemini) to extract the arXiv subcategory
tags from the PDF metadata. By grouping and counting these tags\footnote{I grouped the tags math.AG, math.DG, math.MG and math.SG under ``geometry'';
math.AC, math.CT, math.GR, math.OA, math.QA, math.RA and math.RT under
``algebra''; math.AP, math.CA, math.CV, math.DS, math.FA and math.NA
under ``analysis''; math.AT, math.GN and math.GT under ``topology'';
math.CO under ``combinatorics''; math.NT under ``number theory'';
math.PR and math.ST under ``probability and statistics''; and math.LO
under ``logic and set theory''. arXiv's classification system requires
exactly one category tag per preprint. For arXiv's explanation of
the tags, see https://arxiv.org/category\_taxonomy. } and then comparing the resulting counts, we can obtain a precise
picture of which subfields contributed more or fewer explanation cases
than expected given their prevalence in the corpus. 

The table below shows the percentages of total cases contributed by
each subject area $S$, together with a ``coefficient of explanatory
richness'' statistic obtained by dividing the percentage of $S$-type
explanation cases in the filtered dataset by the percentage of $S$-type
papers in the original corpus. A coefficient of $1$ means that the
subject area punches exactly at its weight, explanatorily speaking.
A value $>1$ or $<1$ means the area is respectively richer or poorer
in explanation cases than expected. 

As the table shows, most subject areas are neither remarkably replete
with nor devoid of appeals to explanation; all coefficients calculated
lay between $0.75$ and $1.35$. Perhaps surprisingly, logic and set
theory earned the largest coefficient, though this may be in part
an artifact of a small sample size (as there were fewer than $100$
papers in this category in the original corpus). Combinatorics is
unusually replete, and geometry is slightly on the rich side. Algebra
and analysis are close to parity with their expected contributions.
Probability and statistics appear to be the most explanatorily impoverished
subfields (though here again the sample was relatively small).\\

\begin{tabularx}{\columnwidth}{>{\centering\arraybackslash}Xccc}
\toprule 
 & \begin{cellvarwidth}[t]
\centering
Percentage of papers \\
in full corpus $=C$
\end{cellvarwidth} & \begin{cellvarwidth}[t]
\centering
Percentage of papers \\
in filtered dataset $=D$
\end{cellvarwidth} & \begin{cellvarwidth}[t]
\centering
``Coefficient of \\
explanatory richness'' $=D/C$
\end{cellvarwidth}\tabularnewline
\midrule
Geometry & $33.5\%$ & $35.6\%$ & $1.06$\tabularnewline
Algebra & $24.2\%$ & $23.9\%$ & $0.99$\tabularnewline
Analysis & $14.3\%$ & $14.3\%$ & $1.00$\tabularnewline
Topology & $10.2\%$ & 9$.0\%$ & $0.88$\tabularnewline
Combinatorics & $5.2\%$ & $6.2\%$ & $1.19$\tabularnewline
Number theory & $3.8\%$ & $3.5\%$ & $0.92$\tabularnewline
Probability and statistics & $3.0\%$ & $2.3\%$ & $0.77$\tabularnewline
Logic and set theory & $1.9\%$ & $2.5\%$ & $1.32$\tabularnewline
\bottomrule
\end{tabularx}\\

\subsection{The prevalence of explanatory concerns}

Philosophers of mathematics have often discussed, and occasionally
debated, whether appeals to explanation are a relatively common and
widespread feature of mathematical practice or a marginal phenomenon.
For instance, \cite{Steiner 1978} asserts that ``{[}m{]}athematicians
routinely distinguish proofs that merely demonstrate from proofs which
explain'' (135), while \cite{Resnik =000026 Kushner 1987} claims
that ``{[}m{]}athematicians rarely describe themselves as explaining''
(151). While subsequent PMP work has produced enough documentary evidence
to fend off extreme versions of skepticism, just how often explanatory
concerns arise in the research literature remains a somewhat open
question.

Gemini's analysis may help clarify this picture. Of the $5000$ papers
in the original corpus, some $735$ contributed at least one example
to the filtered dataset. I estimated above that approximately $20\%$
of these examples involved unambiguous, substantive appeals to explanation,
while another $60\%$ could be interpreted more or less reasonably
as involving explanatory concerns. Assuming (as seems plausible) that
the statistics for the papers roughly follow the statistics for the
examples, this suggests that at least \textasciitilde 150 out of
5000 research papers (around $3\%$) involve clear explanation claims,
while at least \textasciitilde 590 out of 5000 (around $12\%$) feature
cases that are borderline or better. 

Since Gemini isn't an infallible detector, the true numbers are presumably
higher by some amount, though given the model's observed performance
I expect it not to have missed large numbers of high-quality cases.
One should also note that these percentages are averages across a
diverse corpus; as discussed above, talk of explanation is more common
in some subfields than others, so a finer-grained picture would likely
show significantly higher percentages for (e.g.) combinatorics and
lower percentages for (e.g.) probability.

Do numbers in this vicinity mean that mathematicians invoke explanatory
concerns ``routinely'' or ``rarely''? The answer of course depends
on the bar one sets for rarity, and it's unclear what standard would
be appropriate here. Having worked on explanation in a PMP vein for
some time, I myself find these numbers unsurprising. A great deal
of working mathematics consists of prosaic reasoning within inherited
frameworks---sketching variations on standard arguments, recalling
and applying results proved elsewhere, detailing the steps of calculations,
and so on---where epistemic concerns aside from accuracy have little
role to play. 

Perhaps a more interesting question than that of rarity---at least
from the viewpoint of resolving skeptical debates about explanation
as a PMP issue---is whether appeals to explanation are sufficiently
\emph{normal} and \emph{widespread} to constitute a settled practice,
and sufficiently \emph{important} when they arise to influence mathematical
research in meaningful ways. For reasons that philosophers of mathematical
practice have spent decades uncovering, the answer to this question
is clearly yes. I therefore don't share Mizrahi's conviction that
it reflects poorly on PMP priorities if empirical methods suggest
that explicit appeals to explanation are relatively uncommon (\cite{Mizrahi 2020}). 

\subsection{Querying the dataset}

How can philosophers make the best use of Gemini's dataset? One can,
of course, just scroll through its hundreds of pages in search of
interesting examples. But the document can be put to more targeted
uses with additional help from LLM tools.

Suppose that one wants to know, for instance, about the \emph{tradeoffs
}mathematicians sometimes face in weighing explanatorily virtuous
proofs over proofs with other desirable qualities. (Thanks to Marc
Lange for suggesting this topic.) One can investigate by querying
Gemini or another model with the dataset file attached to the prompt.
(E.g.: ``I'm providing you with a set of excerpts from math papers
documenting examples of explanatory proof or discussion of other explanatory
concerns, which you helped produce in an earlier run. I now have some
follow-up questions about the dataset. First: can you find any cases
where mathematicians discuss tradeoffs between explanatory value,
or specific explanatory virtues, and other desirable qualities in
a proof (e.g. simplicity, brevity, elementariness, purity, etc.)?'') 

Gemini produced seven interesting examples when prompted in this way.
One such case: 
\begin{quote}
The reader will probably be relieved to know that, in the cases studied
in this article, no essential use is made of the theory of exterior
differential systems as such. For example, there is no use of the
concepts of polar spaces, regularity, involutivity, characteristic
variety, and so on. ...This avoidance does not significantly lengthen
any of the proofs, so I feel that the savings of not having to introduce
and discuss concepts from exterior differential systems justifies
this strategy. The main disadvantage to the reader is that it does
not explain why the rigidity results could have been anticipated,
making them seem somewhat miraculous. (\cite{Bryant 2001}, 6)
\end{quote}
In another followup I asked Gemini for cases in which explanatory
concerns are expressed without the term `explanation' and its cognates
and obtained a similar number of strong examples. 

Further possibilities abound. One could query one's model of choice
for, say, impure explanations involving combinatorics (``Note that
Proposition 4 shows that, except for an oscillating factor, the creation
rates also behave like $2u$. William Jockusch has shown... how to
use a generating function {[}developed elsewhere{]} to explain this
behavior, by viewing the creation rates as numerical approximations
to a solution of the wave equation'' (\cite{Cohn et al. 2000}, 26)),
or cases where explanatoriness is linked to transparency in the sense
of \cite{D=002019Alessandro 2021} (``There is a lot of information
known on the value of $d\left(A\right)$ for various examples... In
this note, we try to give direct and explicit factorizations for the
corresponding results for $\ell\left(A\right)$... our results are
essentially already known but they shed some light on the meaning
of $\ell\left(A\right)$'' (\cite{Pisier 2000}, 2), or cases which
seem unusual or heterodox in light of extant philosophical work (Gemini
points to \cite{Li 2001}, in which the key insight ``\textsf{is
not just `this proof works', but `this other, more obvious proof
}\textsf{\emph{doesn't }}\textsf{work, and here's the obstacle that
forces us into a more complex strategy'. This justification based
on the failure of simpler methods is a sophisticated form of explanatory
reasoning that is rarely the focus of philosophical case studies}''),
or whatever one's heart desires.

\subsection{Gemini as philosophical assistant}

Since the arrival of advanced LLMs, some philosophers have displayed
interest in the technology's potential to contribute directly to
philosophical discourse. The availability of the present dataset gives
a good opportunity to test Gemini's abilities.

I began by asking Gemini to adjudicate between existing theories of
mathematical explanation in light of the collected cases. (``Which
philosophical theory of explanation seems to best account for the
observed examples, especially those you judge to be clearest, strongest
and most representative? Explain your reasoning.'' The filtered dataset
file was included with this prompt.) 

Gemini declined to endorse a single account, instead arguing that
unificationism, ontic theories and epistemic/pragmatic accounts all
capture important facets of explanation. After a long discussion supported
by many specific cases from the dataset, Gemini summarized its conclusions
as follows:
\begin{enumerate}
\item \textsf{\textbf{Kitcher's Unificationism:}}\textsf{ This theory strongly
resonates with the many examples where mathematicians explain by showing
how specific results fit into broader patterns, derive from general
principles, or connect disparate areas. The drive to reduce the number
of independent phenomena or argument patterns is a powerful explanatory
move in mathematics.}
\item \textsf{\textbf{A General Ontic Account Focused on Revealing Salient
Structures/Properties: }}\textsf{This is crucial for understanding
explanations that pinpoint why a result holds by reference to the
inherent mathematical properties, structures, or mechanisms at play
(e.g., obstructions, characterizing conditions, the nature of a specific
definition).}
\end{enumerate}
\begin{quote}
\textsf{These two, when seen as delivering the }\textsf{\textbf{epistemic
goods}}\textsf{ of increased understanding, clarity, and insight,
seem to provide the most comprehensive account of the observed examples.
The explicit language used by mathematicians strongly supports the
idea that they are seeking explanations that either unify or reveal
the determining underlying structures, precisely because such achievements
lead to a deeper cognitive grasp of the subject matter.}
\end{quote}
I then asked Gemini about problem cases for its preferred accounts.
(``Which cases (or types of case) from the analyses seem to present
the greatest challenge for ontic and unificationist approaches to
explanation? What alternative approaches might better account for
these cases?'') The model gave the following response, culminating
in a tentative endorsement of pluralism or epistemicism (with ontic
or unificationist elements) as the likeliest path forward:
\begin{quote}
\textsf{\textbf{Challenges for Ontic Approaches:}}

\textsf{Ontic accounts generally hold that explanations track objective
dependency relations in the world (or, for mathematics, objective
relations between mathematical facts/structures). The explanans makes
the explanandum happen or be the way it is due to these real relations.}

\textsf{\textbf{Challenging Cases/Types of Cases:}}
\begin{enumerate}
\item \textsf{\textbf{Explanations via ``heuristic interpretation'' or
``intuitive analogy'' that are explicitly not the rigorous proof:}}
\begin{itemize}
\item \textsf{\textbf{Example: }}\textsf{Anshelevich providing a ``heuristic
interpretation'' of the functional Itô formula, linking it to classical
Taylor expansions, explicitly for understanding }\textsf{\emph{before}}\textsf{
the rigorous (and potentially complex) proof. }
\item \textsf{\textbf{Challenge:}}\textsf{ While the heuristic might point
to an underlying structural analogy, the }\textsf{\emph{explanation
itself}}\textsf{ in these instances seems to reside in the cognitive
act of relating the complex to the simple and familiar, even if that
analogy isn't a perfect or rigorous ontic dependency. The value is
in making it ``graspable'', which is more of an epistemic/pragmatic
virtue than a direct reflection of an objective dependency chain in
the formal proof. An ontic account would need to argue that the heuristic
does track some objective (perhaps looser) structural resemblance.}
\end{itemize}
\item \textsf{\textbf{Reproofs explicitly valued for ``simplicity'', ``clarity'',
or being ``less computational'' when the original proof already
established the ontic link:}}
\begin{itemize}
\item \textsf{\textbf{Example:}}\textsf{ Schlosser's ``elementary derivations''
of hypergeometric identities, contrasting with prior ``not entirely
elementary'' proofs. Or the desire for a ``more conceptual'' proof
by Burghelea \& Haller.}
\item \textsf{\textbf{Challenge:}}\textsf{ If the original proof already
correctly identified the mathematical facts from which the theorem
follows (the ontic basis), why is a new, simpler proof more ``explanatory''?
The ontic dependencies might be the same. The added value seems to
be epistemic: the new proof makes those dependencies easier to grasp
or presents them in a more ``natural'' way. An ontic account might
struggle to capture this difference in explanatory value if the objective
dependencies are unchanged.}
\end{itemize}
\item \textsf{\textbf{Explanations where the ``explanatory work'' seems
to be done by a representation or framework rather than directly by
the objects themselves:}}
\begin{itemize}
\item \textsf{\textbf{Example:}}\textsf{ The ``logarithmic point of view''
(El Goul) making a refined result an ``easy consequence''. The framework
itself facilitates understanding. Or Fialowski \& Penkava discussing
how reframing deformation theory in terms of coderivations makes connections
``completely transparent''. }
\item \textsf{Challenge: While the framework undoubtedly reflects ontic
structures, the explanatory power seems to come from the choice of
representation that highlights certain features and relationships.
An ontic account needs to carefully delineate whether the explanation
lies in the objective structures themselves or in the way a particular
framework allows us to cognize those structures.}
\end{itemize}
\end{enumerate}
\textsf{\textbf{Challenges for Unificationist Approaches (Kitcher):}}

\textsf{Unificationism sees explanation as deriving phenomena from
a minimal set of argument patterns. The more phenomena derived from
fewer, more stringent patterns, the better the unification and thus
the explanation.}

\textsf{\textbf{Challenging Cases/Types of Cases:}}
\begin{enumerate}
\item \textsf{\textbf{``Local'' explanations of specific proof steps or
algorithmic choices:}}
\begin{itemize}
\item \textsf{\textbf{Example:}}\textsf{ The explanation in Osofsky's paper
for why a specific pivot in an algorithm is congruent to $1$ modulo
$p\nu$, based on determinant properties and prior row operations.}
\item \textsf{\textbf{Challenge:}}\textsf{ This is a very local explanation
for a specific step. It's not immediately clear how this contributes
to a global unification of a broad range of mathematical facts under
a minimal set of argument patterns. The explanation here seems more
about a local, causal-like chain of deductive reasoning within the
proof.}
\end{itemize}
\item \textsf{\textbf{Explanations that ``demystify'' by providing a single
crucial insight or connection, rather than subsuming under a broad
pattern:}}
\begin{itemize}
\item \textsf{\textbf{Example:}}\textsf{ Garoufalidis \& Levine ``demystifying''
the Gessel/Reutenauer map by making a connection, not necessarily
by showing it's an instance of a widely applicable argument pattern. }
\item \textsf{\textbf{Challenge: }}\textsf{Unificationism values reducing
the number of fundamental patterns. An explanation that provides a
key, novel insight for a specific phenomenon might be highly valued
by mathematicians for its explanatory power, even if it doesn't immediately
contribute to a grand unification scheme by reducing the number of
basic argument patterns used across a wider domain.}
\end{itemize}
\item \textsf{\textbf{Explanations based on identifying why a general approach
fails in a specific case, by pointing to a unique obstructing feature:}}
\begin{itemize}
\item \textsf{\textbf{Example:}}\textsf{ Suciu explaining why an earlier
conjecture about Chen group ranks failed by pointing to an ``essential
component in $R_{1}$'' present in the Pappus arrangement counterexample.}
\item \textsf{\textbf{Challenge: }}\textsf{This type of explanation highlights
disunification or the limits of a pattern. While understanding the
limits of a unifying pattern is important, the explanation for the
failure itself might not fit neatly into deriving the phenomenon from
an accepted argument pattern. It's more about showing why an expected
pattern doesn't apply.}
\end{itemize}
\end{enumerate}
\textsf{\textbf{Alternative Approaches that Might Better Account for
These Cases:}}

\textsf{For the cases challenging purely ontic or unificationist views,
approaches that incorporate }\textsf{\textbf{epistemic and pragmatic
factors}}\textsf{ seem more promising... }{[}here Gemini gives a summary
of explanatory pluralism, understanding-based and pragmatic accounts
of explanation and their utility for handling the above types of case{]}. 

\textsf{\textbf{Synthesis:}}

\textsf{It's unlikely that a single, monolithic theory will account
for all instances of what mathematicians deem ``explanatory''.}
\begin{itemize}
\item \textsf{The drive for unification is undeniably a powerful explanatory
force in mathematics, and Kitcher's account captures this well for
large-scale theoretical connections.}
\item \textsf{The identification of objective structures and dependencies
(ontic accounts) is fundamental to what is being explained and what
does the explaining in a factual sense.}
\item \textsf{However, the clearest challenges arise when the epistemic
value (simplicity, clarity, production of understanding, demystification)
seems to be the primary justification for calling something explanatory,
especially in cases of reproof or when contrasting different conceptual
frameworks.}
\end{itemize}
\textsf{Therefore, a }\textsf{\textbf{pluralistic approach}}\textsf{
that acknowledges different explanatory virtues, or an }\textsf{\textbf{epistemic
approach}}\textsf{ that defines explanation in terms of the cognitive
benefits it provides (often by revealing unifying or ontic structures
in an accessible way), would likely offer the most comprehensive account
for the full range of examples, especially those that prove challenging
for stricter ontic or unificationist views alone. These alternative
or supplementary approaches emphasize the human endeavor of mathematics
and the goal of making its intricate structures comprehensible.}
\end{quote}
How does this discussion rate as a piece of philosophy? Perhaps most
experts would find few of Gemini's points completely novel. But it's
difficult not to admire the clarity with which the model makes its
case, the cogency of the problems it raises, the appropriateness of
its supporting examples and its grasp of subtle features of the dialectic.
It's striking, for instance, that Gemini sketches a plausible line
of defense for the ontic theorist in response to its own challenge
from explanatory heuristics and analogies, suggesting that ``\textsf{{[}a{]}n
ontic account would need to argue that {[}such a{]} heuristic does
track some objective (perhaps looser) structural resemblance}''.
This is reminiscent of recent moves in debates about the epistemology
of unrealistic scientific models, which have seen discussion of whether
and how idealizations might serve as sources of knowledge about objective
dependence relations (\cite{Rice 2016}, \cite{Sullivan =000026 Khalifa 2019},
\cite{D=002019Alessandro 2023}). 

I believe Gemini's discussion contains at least one novel and philosophically
promising suggestion. Talk of \emph{obstructions} (to particular constructions,
proof strategies or the existence of certain objects) is ubiquitous
in contemporary mathematics---indeed, there are well-developed theories
devoted to systematizing knowledge about obstructions in various domains---\footnote{For instance, ``obstruction theory'' in homotopy theory is the study
of obstacles to extending certain mappings defined on simplicial complexes.}but to my knowledge nobody has written at length about the search
for obstructions as a part of explanatory inquiry. Gemini's suggestion
that this is sometimes so strikes me as very plausible. One can see
prospective connections here with, for instance, questions of when
and how impossibility facts can serve as \emph{explanantia} (\cite{Baron =000026 Colyvan 2020}).
And further motivating examples from the mathematical literature aren't
difficult to find. For instance: ``the Brauer-Manin obstruction explains
all violations of the Hasse principle for del Pezzo surfaces over
a number field $k$'' (\cite{Liedtke 2013}, 301); ``{[}a{]}s we
will later see, similar properties {[}to those previously noted{]}
do not hold for generalized quadrangles, not even in the finite case.
This obstruction is the essential reason why several problems concerning
parameters or classification of quadrangles... become very hard to
handle, while they can be handled smoothly for planes'' (\cite{Thas 2011},
275). Following up on the model's observation here strikes me as
an eminently worthwhile PMP project.

Gemini is thus able to perceptively summarize, comment on and apply
the cases in its dataset to PMP questions about explanation. Can it
produce original philosophical ideas of the same level of quality?
I asked it to try. (``Suppose you were asked to come up with an entirely
novel theory of explanation, one able to do better than existing theories
at making sense of all the cases you've reviewed. What might such
a theory look like? Feel free to speculate.'') The model's answer:
\begin{quote}
\textsf{If I were to speculate on an entirely novel theory of mathematical
explanation, aiming to encompass the diverse cases reviewed, I would
try to build something that integrates insights from existing theories
while addressing their limitations. Let's call it }\textsf{\textbf{``Explanatory
Resonance Theory'' (ERT)}}\textsf{.}

\textsf{Here's what ERT might look like:}

\textsf{\textbf{Core Idea of Explanatory Resonance Theory (ERT):}}

\textsf{Mathematical explanation occurs when a target mathematical
statement or phenomenon (the }\textsf{\emph{explanandum}}\textsf{)
is brought into }\textsf{\textbf{``explanatory resonance''}}\textsf{
with a broader, more structured, or more fundamental mathematical
context (the }\textsf{\emph{explanans-framework}}\textsf{). This resonance
is not a single type of connection but a multi-faceted relationship
that enhances our }\textsf{\textbf{cognitive grasp}}\textsf{ and }\textsf{\textbf{structural
integration}}\textsf{ of the explanandum.}

\textsf{\textbf{Key Components and Features of ERT:}}
\begin{enumerate}
\item \textsf{\textbf{Multi-faceted Resonance:}}\textsf{ ``Resonance''
isn't a monolithic concept. It can manifest in several ways, acknowledging
the pluralistic nature of what mathematicians value as explanatory:
}{[}here Gemini lists and describes ``structural'', ``unifying'',
``simplifying/clarifying'', ``analogical/metaphorical'' and ``generative''
varieties of resonance, attempting to incorporate the core ideas of
the views discussed previously and bringing up many of its earlier
examples{]}.
\item \textsf{\textbf{The Explanans-Framework:}}\textsf{ This is not just
a single proposition but a richer context. It could be:}
\begin{itemize}
\item \textsf{A well-established theory (e.g., Galois theory, representation
theory, category theory).}
\item \textsf{A specific set of definitions and their consequences.}
\item \textsf{A powerful visualization or diagrammatic system.}
\item \textsf{A newly introduced concept or structure that re-organizes
existing knowledge.}
\end{itemize}
\item \textsf{\textbf{Cognitive Grasp and Structural Integration (The Goal):}}\textsf{
The ultimate aim of explanation under ERT is to enhance the mathematician's
cognitive grasp of the explanandum. This means:}
\begin{itemize}
\item \textsf{\textbf{Understanding Why:}}\textsf{ Not just }\textsf{\emph{that}}\textsf{
it's true, but }\textsf{\emph{why}}\textsf{ it must be true given
the explanans-framework.}
\item \textsf{\textbf{Seeing Connections:}}\textsf{ Understanding how the
explanandum relates to other mathematical ideas and structures.}
\item \textsf{\textbf{Increased Predictive/Manipulative Power (Potentially):}}\textsf{
A good explanation might suggest generalizations, new questions, or
more efficient ways to work with the concepts.}
\item \textsf{\textbf{Reduced Surprise/Mystery:}}\textsf{ Phenomena that
initially seem ``surprising'' or ``mysterious'' become less so
when they resonate with a suitable framework.}
\end{itemize}
\item \textsf{\textbf{Dynamic and Contextual:}}\textsf{ What counts as a
good explanans-framework or a strong resonance can be:}
\begin{itemize}
\item \textsf{\textbf{Field-dependent:}}\textsf{ Different mathematical
subfields might value different types of resonance. Bijective proofs
are highly valued in combinatorics for their structural resonance.}
\item \textsf{\textbf{Historically contingent:}}\textsf{ What was once a
complex explanandum can become part of a familiar explanans-framework
for new explananda.}
\item \textsf{\textbf{Audience-dependent (to some extent):}}\textsf{ While
striving for objectivity, the effectiveness of an explanation in producing
``cognitive grasp'' can have a pragmatic element related to the
audience's background knowledge.}
\end{itemize}
\end{enumerate}
\textsf{\textbf{How ERT Might Do Better Than Existing Theories:}}

{[}Here Gemini rehearses the virtues of ERT, e.g. that it accommodates
pluralism, integrates the ontic and epistemic, makes sense of reproof
phenomena, handles heuristic explanations, accounts for explanations
of failure, and emphasizes the role of frameworks and representations
in explanation.{]}

\textsf{\textbf{Potential Challenges for ERT:}}
\begin{itemize}
\item \textsf{\textbf{Defining ``Resonance'' Precisely:}}\textsf{ The
term is evocative but would need further philosophical work to be
made more precise for its different facets. How do we measure or compare
different types/strengths of resonance?}
\item \textsf{\textbf{Objectivity vs. Subjectivity:}}\textsf{ While aiming
to ground explanation in objective mathematical relationships, the
emphasis on ``cognitive grasp'' might introduce an element of subjectivity.
ERT would need to balance this.}
\item \textsf{\textbf{Demarcation:}}\textsf{ When does a connection become
strong enough to count as ``explanatory resonance'' versus mere
logical entailment or a weak analogy?}
\end{itemize}
\textsf{In essence, Explanatory Resonance Theory would propose that
mathematical explanation is the process of making a mathematical idea
``vibrate sympathetically'' with a larger, more harmonious, and
well-understood part of the mathematical orchestra. The quality of
the explanation depends on the richness, depth, and clarity of this
induced resonance.}
\end{quote}
While Gemini's sophisticated grasp of the issues remains on display
here, its efforts to stake out a novel view strike me as less impressive.
As the model itself recognizes, the core idea of ERT is a metaphor
whose intuitive content is far from clear (and which is far removed
from the way mathematicians tend to describe their explanatory practices).
Indeed, Gemini seems to use ``resonance'' as a catchall for any
sort of relationship between mathematical theories, facts, frameworks
or phenomena which might in some cases be associated with explanatory
understanding. In general, It's difficult to see what such a broad
and ill-defined construct might add to existing versions of epistemicism
or pluralism.

To be sure, the above was a single answer to a single prompt; one
might very well wonder whether the model could do better with further
prompting and more space to expand on its views. My experience posing
additional critical and clarificatory questions was that Gemini continued
to propose lists of (fluent and reasonable but standard) ideas which
might be helpful for developing ERT, without committing to much of
anything or producing any recognizably new insights. This result aligns
with the widespread perception that frontier language models as yet
lack the cognitive skills to make substantive intellectual discoveries
(cf. \cite{Altman 2025}; Altman, CEO of OpenAI, opines that ``2026
will likely see the arrival of systems that can figure out novel insights''). 

\section{Discussion}

This paper has shown something of what current AI systems can offer
to PMP researchers. First and most notably, they've laid a path for
a novel paradigm in corpus analysis. Improvements in content length,
chain-of-thought reasoning and hallucination control enable current
models to process large volumes of mathematical text, and to report
comprehensively, intelligently and accurately about findings of interest.
These methods can help practice-oriented philosophers build sizable
datasets of examples for use in future work. The explanation dataset
created here, consisting of well over $1,000$ annotated references
from over $700$ individual papers, illustrates the scope of what's
possible. 

Current models can moreover help researchers automate many labor-intensive
processes, by authoring scripts to interact with the models' own
APIs and perform other useful tasks. As a programming novice, Gemini's
advanced code-writing capabilities were indispensable to many aspects
of the present study, from downloading the arXiv corpus to queuing
up hundreds of automated queries to extracting information from PDF
metadata. The availability of these tools can go far in empowering
humanities researchers to execute large projects with technical elements.

AI tools can also help researchers engage philosophically with the
contents of large datasets. Models like Gemini have proven proficient
in understanding and reasoning about concepts of interest in specialist
debates. The present study also shows that they can apply these concepts
skillfully to repositories of examples. At minimum, given appropriate
prompting and context, it's now reasonable to expect sophisticated
feedback from frontier models on how a complex dataset might support
or challenge a philosophical theory or otherwise inform a debate.
With luck one might even get good ideas about novel research directions,
as in the case of Gemini's observation about obstructions and explanation.

This experiment represents a first step along the path of AI-assisted
PMP research. Its results could undoubtedly be improved upon or added
to in a variety of ways: through better prompting and filtering, deeper
exploration of the dataset, more extensive investigation of the model's
analytical capabilities, and so on. I present this study as a proof
of concept in hopes that others will be encouraged to learn from my
experiences and push the limits further. There are of course countless
topics besides explanation in PMP and beyond to which these methods
could be fruitfully applied. 

What concrete steps might philosophers wish to take next? Here are
a few suggestions. First, it would be worthwhile to repeat some or
all of this work with other leading models, both to obtain alternative
perspectives on the subject matter and to gain a better understanding
of the strengths and weakness of existing tools. At the time of writing,
Claude Opus 4, various of OpenAI's models (o3, GPT-4.1 and -4.5)
and DeepSeek-R1 offer comparable performance to Gemini 2.5 (although
only GPT 4.1 matches its one million token context length). Second,
more experimentation with prompts would likely repay the effort. It
would be useful to know, for instance, to what extent explicit examples
of desired and undesired classification behavior improves performance
(and how many examples are optimal), what sorts of system prompts
or other high-level instructions do the most good, and how much insight
the model gains from long-context scholarly sources like the \emph{Stanford
Encyclopedia }excerpt used here. Third, granting models internet access
for search, code execution permission to help with complex tasks and
other scaffolding might prove helpful with certain research tasks.
(I chose not to enable Google search functionality for my interactions
with Gemini, and code execution is currently unavailable for the 2.5-series
models.) Fourth, some may wish to repeat something like the present
study with an even larger corpus. My 5000-paper sample is still a
fraction of the full (freely available) arXiv mathematics repository
of some 80,000 texts. There are likely lessons to be learned from
the greater chronological and subject-matter diversity present in
a larger sample. Fifth, one could instead work exclusively with published
texts selected using some criteria of interest. This approach would
for instance allow for more meaningful comparative work, e.g. examining
longitudinal trends in leading journals over many decades.

Finally, what technical developments should philosophers look to for
improving the prospects of this type of research? Contexts windows
still larger than the current state of the art would be extremely
useful. (Indeed, they'll be a practical necessity for studies much
larger than this one: if a 5000-paper corpus generates over 500 pages
of analysis text, a corpus five or ten times the size will produce
a larger dataset than current models can digest in one sitting.) A
model better able to assess its own confidence level about the quality
of reported examples could help create datasets with built-in filtering.
Further down the line, perhaps, AI research platforms with larger
suites of integrated tools could more seamlessly handle parts of the
project workflow that remain human-labor-intensive (e.g., autonomously
accessing a specified digital corpus without the need for user-side
downloading, scripting, batching and the like). Likewise, enhanced
agentic planning and tool use capabilities would allow models to design
and execute more complex research projects.

How AI systems will develop in the coming years, and how the availability
of these systems will affect philosophy as a whole, are difficult
matters to predict. Even if the frontier of expert-level creative
intelligence remains far off, however, practice-oriented researchers
have been handed a valuable tool in today's crop of LLMs. I look
forward to the progress they'll make possible.


\begin{thebibliography}{D’Alessandro \& Lehet 2024}
\bibitem[Alekseev et al. 2000]{Alekseev et al. 2000}Alekseev, Anton,
Eckhard Meinrenken and Chris Woodward. 2000. ``Formulas of Verlinde
type for non-simply connected groups.'' arXiv:math/0005047.

\bibitem[Allcock et al. 2000]{Allcock et al. 2000}Allcock, Daniel,
James A. Carlson and Domingo Toledo. 2000. ``The complex hyperbolic
geometry of the moduli space of cubic surfaces.'' arXiv:math/0007048.

\bibitem[Altman 2025]{Altman 2025}Altman, Sam. 2025. ``The gentle
singularity.'' \emph{Sam Altman}, June 10 2025 (accessed June 13
2025), URL = https://blog.samaltman.com/the-gentle-singularity.

\bibitem[Baron \& Colyvan 2020]{Baron =000026 Colyvan 2020}Baron,
Sam and Mark Colyvan. 2020. ``Explanation impossible.'' \emph{Philosophical
Studies} 178, 559--576.

\bibitem[Bessis et al. 2000]{Bessis et al. 2000}Bessis, David, Francois
Digne and Jean Michel. 2000. ``Springer theory in braid groups and
the Birman-Ko-Lee monoid.'' arXiv:math/0010254.

\bibitem[Bressoud 2000]{Bressoud 2000}Bressoud, David M. 2000. ``Three
alternating sign matrix identities in search of bijective proofs.''
arXiv:math/0007114.

\bibitem[Bryant 2001]{Bryant 2001}Bryant, Robert L. 2001. ``Rigidity
and quasi-rigidity of extremal cycles in Hermitian symmetric spaces.''
arXiv:math/0006186.

\selectlanguage{english}%
\bibitem[Chartrand 2022]{Chartrand 2022}Chartrand, Louis. 2022. ``Modeling
and corpus methods in experimental philosophy.'' \emph{Philosophy
Compass} 17 (6), 2022.

\selectlanguage{american}%
\bibitem[Cherednik 2000]{Cherednik 2000}\foreignlanguage{english}{Cherednik,
Ivan. 2000. ``One-dimensional double Hecke algebras and Gaussians.''
arXiv:math/0003017.}

\bibitem[Cohn et al. 2000]{Cohn et al. 2000}\foreignlanguage{english}{Cohn,
Henry, Noam Elkies and James Propp. 2000. ``Local statistics for
random domino tilings of the Aztec diamond.'' arXiv:math/0008243.}

\bibitem[D’Alessandro 2020]{D=002019Alessandro 2020}\foreignlanguage{english}{D'Alessandro,
William. 2020. ``Mathematical explanation beyond explanatory proof.''
\emph{British Journal for the Philosophy of Science} 71, 581--603.}

\bibitem[D’Alessandro 2021]{D=002019Alessandro 2021}\foreignlanguage{english}{D'Alessandro,
William. 2021. ``Proving quadratic reciprocity: Explanation, disagreement,
transparency and depth.'' \emph{Synthese} 198, 8621--8664.}

\bibitem[D’Alessandro 2023]{D=002019Alessandro 2023}\foreignlanguage{english}{D'Alessandro,
William. 2023. ``Unrealistic models in mathematics.'' \emph{Philosophers'
Imprint} 23, DOI: 10.3998/phimp.1712.}

\bibitem[D’Alessandro 2025]{D=002019Alessandro 2025}D'Alessandro,
William. 2025. ``Toward a methodology for the philosophy of mathematical
practice.'' \emph{Philosophy of Science}, DOI: XXXX.

\bibitem[D’Alessandro \& Lehet 2024]{D=002019Alessandro =000026 Lehet 2024}D'Alessandro,
William and Ellen Lehet. 2024. ``A noetic account of explanation
in mathematics.'' \emph{Philosophical Quarterly}, DOI: 10.1093/pq/pqae137.

\bibitem[Dolgachev \& Keum 2001]{Dolgachev =000026 Keum 2001}Dolgachev,
Igor V. and Jonghae Keum. 2001. ``Birational automorphisms of quartic
Hessian surfaces.'' arXiv:math/0103186.

\bibitem[Fill \& Janson 2000]{Fill =000026 Janson 2000}Fill, James
Allen and Svante Janson. 2000. ``Smoothness and decay properties
of the limiting Quicksort density function.'' arXiv:math/0005235.

\bibitem[Gannon 2002]{Gannon 2002}Gannon, Terry. 2002. ``Modular
data: The algebraic combinatorics of conformal field theory.'' arXiv:math/0103044.

\bibitem[Guillemin \& Zara 2000]{Guillemin =000026 Zara 2000}Guillemin,
Victor and Catalin Zara. 2000. ``Combinatorial formulas for products
of Thom classes.'' arXiv:math/0007166.

\selectlanguage{english}%
\bibitem[Inglis \& Aberdein 2015]{Inglis =000026 Aberdein 2015}Inglis,
Matthew and Andrew Aberdein. 2015. ``Beauty is not simplicity: An
analysis of mathematicians' proof appraisals.'' \emph{Philosophia
Mathematica }23, 87-109.

\bibitem[Inglis \& Aberdein 2016]{Inglis =000026 Aberdein 2016}Inglis,
Matthew and Andrew Aberdein. 2016. ``Diversity in proof appraisal.''
In Brendan Larvor (ed.), \emph{Mathematical Cultures}, Birkhäuser:
Switzerland, 163--179.

\selectlanguage{american}%
\bibitem[Johansen et al. 2018]{Johansen et al. 2018}\foreignlanguage{english}{Johansen,
Mikkel Willum, Morten Misfeldt and Josefine Lomholt Pallavicini. 2018.
``A typology of} \foreignlanguage{english}{mathematical diagrams.''
In Peter Chapman, Gem Stapleton, Amirouche Moktefi, Sarah Perez-Kriz}
\foreignlanguage{english}{\& Francesco Bellucci (eds.), \emph{Diagrammatic
Representation and Inference: Proceedings of the 10th International
Conference, Diagrams 2018}, Springer. DOI: 10.1007/978-3-319-91376-6\_13.}

\bibitem[Kaloshin 2000]{Kaloshin 2000}\foreignlanguage{english}{Kaloshin,
Vadim. 2000. ``A geometric proof of the existence of Whitney stratifications.''
arXiv:math/0010144.}

\selectlanguage{english}%
\bibitem[Lange 2017]{Lange 2017}Lange, Marc. 2017. ``Mathematical
explanations that are not proofs.'' \emph{Erkenntnis}, DOI: 10.1007/s10670-017-9941-z.

\bibitem[Lehet 2021]{Lehet 2021}Lehet, Ellen. 2021. ``Mathematical
explanation in practice.'' \emph{Axiomathes} 31, 553--574.

\selectlanguage{american}%
\bibitem[Li 2001]{Li 2001}\foreignlanguage{english}{Li, Jun. 2001.
``A degeneration of stable morphisms and relative stable morphisms.''
arXiv:math/0009097.}

\bibitem[Liedtke 2013]{Liedtke 2013}\foreignlanguage{english}{Liedtke,
Christian. 2013. ``Algebraic surfaces in positive characteristic.''
In Fedor Bogomolov, Brendan Hassett and Yuri Tschinkel (eds.), \emph{Birational
Geometry, Rational Curves, and Arithmetic}, Springer: New York, 229--292.}

\bibitem[Loeb \& Talvila 2001]{Loeb =000026 Talvila 2001}\foreignlanguage{english}{Loeb,
Peter A. and Erik Talvila. 2001. ``Covering theorems and Lebesgue
integration.'' arXiv:math/0101014.}

\bibitem[Mancosu et al. 2023]{Mancosu et al. 2023}\foreignlanguage{english}{Mancosu,
Paolo, Francesca Poggiolesi and Christopher Pincock. 2023. ``Mathematical
explanation.'' In Edward N. Zalta and Uri Nodelman (eds.), \emph{The
Stanford Encyclopedia of Philosophy }(Fall 2023 Edition), URL = <https://plato.stanford.edu/archives/fall2023/entries/mathematics-explanation/>.}

\bibitem[Massey 2001]{Massey 2001}\foreignlanguage{english}{Massey,
David B. 2001. ``A little microlocal Morse theory.'' arXiv:math/0006185.}

\bibitem[McDuff 2004]{McDuff 2004}\foreignlanguage{english}{McDuff,
Dusa. 2004. ``Symplectomorphism groups and almost complex structures.''
arXiv:math/0010274.}

\bibitem[Mejía Ramos et al. 2019]{Mej=0000EDa Ramos et al. 2019}\foreignlanguage{english}{Mejía
Ramos, Juan Pablo, Lara Alcock, Kristen Lew, Paolo Rago, Chris Sangwin
and Matthew Inglis. 2019. ``Using corpus linguistics to investigate
mathematical explanation.'' In Eugen Fischer and Mark Curtis (eds.),
\emph{Methodological Advances in Experimental Philosophy}, London:
Bloomsbury, 239--263.}

\selectlanguage{english}%
\bibitem[Mejía Ramos et al. 2021]{Mej=0000EDa Ramos et al. 2021}Mejía
Ramos, Juan Pablo, Tanya Evans, Colin Rittberg and Matthew Inglis.
2021. ``Mathematicians' assessments of the explanatory value of
proofs.'' \emph{Axiomathes} 31, 575--599.

\selectlanguage{american}%
\bibitem[Milson 2000]{Milson 2000}\foreignlanguage{english}{Milson,
Robert. 2000. ``Composition sum identities related to the distribution
of coordinate values in a discrete simplex.'' arXiv:math/0003126.}

\selectlanguage{english}%
\bibitem[Mizrahi 2020]{Mizrahi 2020}Mizrahi, Moti. 2020. ``Proof,
explanation, and justification in mathematical practice.'' \emph{Journal
for General Philosophy of Science} 51, 551-568.

\selectlanguage{american}%
\bibitem[Pease et al. 2019]{Pease et al. 2019}\foreignlanguage{english}{Pease,
Alison, Andrew Aberdein and Ursula Martin. 2019. ``Explanation in
mathematical conversations: An empirical investigation.'' \emph{Philosophical
Transactions of the Royal Society A }377: 20180159. DOI: 0.1098/rsta.2018.0159.}

\bibitem[Pisier 2000]{Pisier 2000}\foreignlanguage{english}{Pisier,
Gilles. 2000. ``Remarks on the similarity degree of an operator algebra.''
arXiv:math/0009052.}

\bibitem[Ramakrishnan 2001]{Ramakrishnan 2001}\foreignlanguage{english}{Ramakrishnan,
Dinakar. 2001. ``Modularity of solvable Artin representations of
$\mathrm{GO}\left(4\right)$-type.'' arXiv:math/0102231.}

\selectlanguage{english}%
\bibitem[Resnik \& Kushner 1987]{Resnik =000026 Kushner 1987}Resnik,
Michael D. and David Kushner. 1987. ``Explanation, independence and
realism in mathematics.'' \emph{British Journal for the Philosophy
of Science }38, 141-158.

\selectlanguage{american}%
\bibitem[Rice 2016]{Rice 2016}\foreignlanguage{english}{Rice, Collin.
2016. ``Factive scientific understanding without accurate representation.''
\emph{Biology and Philosophy} 31, 81--102. }

\selectlanguage{english}%
\bibitem[Rittberg 2019]{Rittberg 2019}Rittberg, Colin Jakob. 2019.
``On the contemporary practice of philosophy of mathematics.'' \emph{Acta
Baltica Historiae et Philosophiae Scientiarum} 7, 5--26.

\selectlanguage{american}%
\bibitem[Sa et al. 2024]{Sa et al. 2024}\foreignlanguage{english}{Sa,
Rentuya, Lara Alcock, Matthew Inglis and Fenner Stanley Tanswell.
2024. ``Do mathematicians agree about mathematical beauty?'' \emph{Review
of Philosophy and Psychology} 15, 299--325.}

\selectlanguage{english}%
\bibitem[Steiner 1978]{Steiner 1978}Steiner, Mark. 1978. ``Mathematical
explanation.'' \emph{Philosophical Studies }34, 135-151.

\selectlanguage{american}%
\bibitem[Sullivan \& Khalifa 2019]{Sullivan =000026 Khalifa 2019}\foreignlanguage{english}{Sullivan,
Emily and Kareem Khalifa. 2019. ``Idealizations and understanding:
Much ado about nothing?'' \emph{Australasian Journal of Philosophy}
97, 673--689.}

\bibitem[Tao 2000]{Tao 2000}\foreignlanguage{english}{Tao, Terence.
2000. ``Global regularity of wave maps I. Small critical Sobolev
norm in high dimension.'' arXiv:math/0010068.}

\bibitem[Thas 2011]{Thas 2011}\foreignlanguage{english}{Thas, Koen.
2011. ``Order in building theory.'' In Robin Chapman (ed.), \emph{Surveys
in Combinatorics 2011}, Cambridge: Cambridge University Press, 235--332.}

\bibitem[van der Kallen 2000]{van der Kallen 2000}\foreignlanguage{english}{van
der Kallen, Wilberd. 2000. ``From Mennicke symbols to Euler class
groups.'' arXiv:math/0010226.}

\end{thebibliography}
\end{document}